\documentclass[11pt]{amsart}

\usepackage{amsmath}
\usepackage{amssymb}
\usepackage{mathrsfs}

\theoremstyle{definition}

\theoremstyle{remark}

\newcount\skewfactor
\def\mathunderaccent#1#2 {\let\theaccent#1\skewfactor#2
\mathpalette\putaccentunder}
\def\putaccentunder#1#2{\oalign{$#1#2$\crcr\hidewidth
\vbox to.2ex{\hbox{$#1\skew\skewfactor\theaccent{}$}\vss}\hidewidth}}

\begin{document}

\title {Model theory }
\author {Saharon Shelah}
\address{Einstein Institute of Mathematics\\
Edmond J. Safra Campus, Givat Ram\\
The Hebrew University of Jerusalem\\
Jerusalem, 91904, Israel\\
 and \\
 Department of Mathematics\\
 Hill Center - Busch Campus \\
 Rutgers, The State University of New Jersey \\
 110 Frelinghuysen Road \\
 Piscataway, NJ 08854-8019 USA}
\email{shelah@math.huji.ac.il}
\urladdr{http://shelah.logic.at}
\thanks{This research was supported by the United States-Israel
Binational Science Foundation/ the Israel basic research foundation/
the Germany-Israel Binational Foundation.
Paper E70. 
We thank  J. Baldwin, G Cherlin, U.Hrushovski,
M.Malliaris and J. Vaananen for helpful comments
}

\date{2011-12-19}

\maketitle
\numberwithin{equation}{section}

\section{Why am I  interested in model theory
(a branch of mathematical logic)?}

In elementary school, mathematics looked (to me)
like just a computational skill---how to multiply,
how to find formulas for areas of squares, rectangles,
triangles etc---and the natural sciences
looked more attractive. Then, entering the ninth
grade, Euclidean geometry captured my heart:
from the  bare bones of assumptions a magnificent
structure is built; an intellectual
endeavour in which it is enough to be right.

Undergraduate mathematics was impressive
for me,
but 
algebra considerably more so than analysis.
Reading  Galois theory, understanding equations in
 general fields, was a gem.
Finding  order in what looks like a
chaos, not grinding water but finding
natural definitions and hard theorems;
generality, being able to say something
from very few assumptions,
was impressive. From this perspective
mathematical logic was the most general direction,
so  I  took the trouble to
do my M.Sc. thesis in mathematical logic;
the thesis happened to be on the model theory
of infinitary logics.

Model theory seemed
the epitome of what I was looking for: rather
than investigating a specific class like ``the class
of fields", the ``class of rings with no zero divisors"
or whatever, we have a class of structures,
called here models.
For this to be meaningful, we have to restrict
somewhat the class, first by saying they are
all of the same ``kind", i.e. have the same
function symbols (for rings: addition, multiplication;
also the so-called ``individual constants"
$ 0$  and $ 1 $, we may have so called predicates,
i.e. symbols for relations, but we shall ignore
that  
point; this information is called the
vocabulary). We have to further restrict the
classes we consider, and the classical
choice in model theory is to restrict to the so called e.c.  i.e. elementary classes,
explained below.

Naturally model theorists start from the bottom:
Consider
$ K $, an e.c.
 (elementary class), i.e.
the class of models of a
first order theory $ T $ 
as explained below. The class
$ K $ (i.e.
$ T $) is called categorical in the infinite
cardinal $ \lambda $ if
it has a unique model up to isomorphism	
of cardinality (= number of elements)
$ \lambda $.
{\L}o\'s  
conjectured  that if  an e.c.
$ K $ with countable vocabulary
is categorical in one uncountable cardinal
then this holds for every uncountable
cardinal.
 After more than a decade, Morley
proved this, and when I started my PhD studies
I thought it was wonderful (and still think so).

The point of view explained above
naturally  leads to the  classification
program. The basic thesis of the classification program is that
reasonable families of classes
of mathematical structures should have natural dividing
lines.
Here a dividing line means
a partition into low, analyzable, tame classes
on the one hand, and high, complicated, wild classes
on the other.
These partitions  will generate a tameness
hierarchy. For each such partition, if the class is on the tame
side one should
have useful structural
analyses applying to all structures in
the class, while if the class is on the wild side one
should have strong evidence of chaotic
behavior (set theoretic complexity).
These results should be complementary, proving that the dividing
lines are not merely sufficient conditions for
being low complexity, or sufficient conditions for being high
complexity.
This calls for relevant test questions; we expect
not to start with a picture of the
meaning of ``analyzable"  and look for a general
context, as this usually does not  provide  evidence for
this being a dividing line.
Of course, although it is hard to refute this thesis
(as you may have chosen the wrong
test questions; in fact this is the nature
of a  thesis), it may lead us to fruitful or
unfruitful directions. The thesis implies
the natural expectation that a success
in developing a worthwhile theory
will lead us also to applications in other parts of
mathematics, but for me this was
neither a prime motivation nor a major
test, just a welcome and
not surprising (in principle)
side benefit
and a ``proof for the uninitiated",
 so we shall not deal with such  %
 important applications.  %

We still have to define what  an e.c. (elementary
class) is. It is
a ``class
of structures satisfying a fixed first
order theory $ T $".
For our purpose, this can be
explained as follows:
given a structure $ M $, we
consider subsets
of $ M $,  sets of pairs
of elements of $ M $, and more generally sets of
 $ n $-tuples of elements of $ M $, 
{\it which are reasonably definable}.  %
By this we mean the following:  %
start with  
structures satisfying an equation (or
another atomic  formula if we have also relation
symbols). Those we call the atomic relations.
 But we may also look at the set of parameters
for which an equation has a solution. More generally,
the set of first order definable relations on
$ M $ is the closure of the atomic ones, under
union (i.e. demanding at least one 
of two conditions, logically ``or")
and intersection (i.e. ``and" ),  under complement
and lastly
we close under projections,
which means ``there is $ x $
such that\dots";
but we do \underline{not} use
\text{``}there is a set of elements\text{''} 
or even \text{``}there is a finite sequence 
of elements\text{''}.
  The way we define such a
set is called a first order
formula, denoted by
$ \varphi (x_0, \dots , x_{n-1}) $.
If $ n=0 $ this will be just true or false
in the structure
and such formulas are called sentences.
The (complete first order) theory
$ Th(M)$  of $ M $
is the set of (first order) sentences it satisfies.
An e.c.(=elementary class) 
 is the class $ Mod_T $
of models of $ T $, that is the
structures $M$ (of the relevant kind, vocabulary)
such that $ Th(M)=T $.
Naturally, N is an elementary  %
extension of
$ M $ (and $ M $ is an elementary submodel of 
$ N $) when for any of those definition,
on finite sequences from %
the smaller model
they agree.
There are many natural classes which are of this form,
 ranging from Abelian
groups and algebraically closed
fields, through random graphs to Peano
Arithmetic, Set Theory, and the like.

A reader may well say
that this setting
 is  too general, that it is nice to deal with
``everything",
but if what we can say is
``nothing",
null or just
dull, then it is not interesting.
 However, this is not the case.
The classification program
 has been successfully done
for the partition to stable/unstable  and
further subdivisions have been established on the tame side
 for the family of  elementary classes.
 Critical dividing lines
for
the taxonomy
involve the
behavior of the Boolean algebras of
parametrically first order definable sets and
relations, i.e.:
$ \varphi (M, \bar{ a } ) =^{\rm df }
\{ \bar{ b } : M $ satisfies
$ \varphi ( \bar{ b } , \bar{ a } )\}  $. 
E.g.  %
$ T = Th(M) $. i.e. $ K=Mod_T $ 
 is unstable iff some first order formula
$ \varphi ( x, y ) $ linearly orders some infinite
set of elements (not necessarily definable itself!)
in some model from %
 $ K $, or similarly for a set of pairs
or, more generally, a set of $ n $-tuples.
A prominent test question involves the number
of models from %
 $ K $ %
 up to isomorphism
of cardinality (= number of elements) $ \lambda $,
called $ I( \lambda , K ) $.  %
The promised ``analyzable" classes include  in this case
 notions of independence
and of dimension (mainly as in the dimension of
a vector space),
and (first order definable) groups and fields
appearing \text{``}out of nowhere\text{''}.

 Clearly having many non-isomorphic
models is a kind of
\text{``}set-theoretic witness for complexity\text{''}
but certainly not a unique one.

Of course what looks like a small corner,
 a family of well understood classes
from the present point of view, looks like a huge cosmos full
of deep mysteries from another point of view,
and  some of these mysteries have resulted in great achievements.

\section{What are, in my opinion,
the most challenging problems in
model theory?}

We may think that the restrictions to
elementary classes is too strong,
but then what takes the place
of the first order definable
parametrized sets? Naturally, at
least a posteriori, we may generalize
this concept but we may wonder
can we really dispense with it?
See \S(2.1).

Dually, we may feel
 that as successful as the dividing line
stable/unstable (and finer divisions "below" that) has been,
not all unstable classes are completely wild,
(and what constitutes being complicated, 
un-analysable depends  on your yard-stick).
Moreover, though many elementary classes 
are stable, many mathematically useful ones are not.

In fact there are provably just two 
\text{``}reasons\text{''}  for being unstable.
The two \text{``}minimal\text{''} 
unstable elementary classes correspond to
the theory of dense linear orders and 
the theory of random graphs.
Much attention has been given
on the one hand to so-called {\it
simple theories}
which include the random graphs and also "pseudo-finite
fields" (see \S(2.1)), and on the other hand
to the {\it dependent theories}, which include the theories
of 
dense linear order, the  %
real field, the $p$-adics, and many fields
of power series (see \S(2.3)).

It is tempting to
 look for a ``maximal (somewhat) tame 
family of
elementary  
classes''.
A natural candidate for this is the following.

We may look for an extreme
 condition of the form of unstability; such a condition
is ``$K  
$ is straight maximal", which means that for some
formula  $ \varphi ( x, y ) $
(or  $ \varphi ( \bar{ x },\bar{ y } ) $),
 for every $ n $ and
non-empty subset  $ {\mathscr F} $
of
${\mathscr F} _n =^{\rm df }  \{ f: f $
a function from
$ \{ 0, \dots , n-1 \} $ to $ \{ 0,1\} \}  $
we can find a model 
$ M  \in K $ 
and $ b_0, \dots , b_{n-1} \in M $
such that:
if $ f \in {\mathscr F} _n $ then
there is $ a \in M $ such that
\text{``}M satisfies $ \varphi (a, b_i ) $
iff $ f(i)=1 $\text{''}
\underline{iff} $ f \in {\mathscr F} $.
Does this really define
an (interesting)
 dividing line?
I am sure it does, and
that
we can say
many things about it;
unfortunately I have neither idea what
those things are, nor of any natural test problem; so we
shall look  instead at problems which have been 
somewhat  
clarified.

\subsection{Non-elementary classes}

 We may note that the family of elementary classes
is a quite restricted family of classes, and many  mathematically natural classes
cannot be described by first order  conditions.
For example ``locally finite'' structures,  such as groups in which
where every finitely generated
subgroup is finite, (or, similarly, is solvable or
the like), or structures satisfying various chain
conditions are not elementarily axiomatizable.

So a fundamental question is
``have a generalization of the
existing stability theory
for a really wide family of  classes,
 where the basic  methods of
e.c. 
completely fail
(in particular, 
nothing like 
the family of parametrically
first order definable sets%
),
but we still have the same test question''.

A good candidate for this broader context is the family of
{\it aec (abstract elementary classes)},
$ {\mathfrak k}   = (K, \le_  {\mathfrak k}   ) $
where $ K $ is a class of models of a fixed vocabulary,
$ \le _ {\mathfrak k}   $ is a partial order on the class
refining the sub-model relation, and
satisfying the obvious properties of e.c.-s
(which means, closure of $ K $ and
$ \le _ {\mathfrak k}   $ are closed under
isomorphism,
any directed system there is
a $ \le _ {\mathfrak k}   $-lub,
every member can be approximated by
$ \le _ {\mathfrak k}   $-submodels of
cardinality bounded by some
$ \chi = LST( {\mathfrak k}   ) $
and "if $ M_ 1 \subseteq M_ 2 $
are $  \le _{\mathfrak k}   $-sub-models
of some $ N $ then $ M_1
\le_ {\mathfrak k}     M_ 2 $.
 (For example, 
those defined by infinitary logics like
the so-called
$   {L}   _{\lambda^+, {\aleph_0}}  $
where we allow conjunctions of $ \lambda $ formulas
but not quantification over infinitely many
variables)). %

Here  a natural  test question is the
large scale, asymptotic behaviour of
$  I(\lambda, K) $, the  
number of 
 models in $ K $ of cardinality
$ \lambda $ up to isomorphism.
Our dream is to prove the main
gap conjecture in this case (see
\S(2.4) below).

The simplest case is the %
 categoricity conjecture:
having
a unique model up to isomorphism 
for every large enough cardinal or
failure of this 
in every large enough cardinal; in-spite of some
advances we still do not know even this,
but there are indications
that a positive theory along these lines exists.

\subsection{Unstable elementary classes- friends of random graphs}

Being a {\it simple} e.c. can be defined similarly
to stability, by ``no first order formula
$ \varphi ( \bar{ x } , \bar{ y } ) $
represents a tree'' (rather than a linear order, in the case of stability)''.
For simple  e.c.-s. 
we know much on analogs of the stable case,
as well as something on non-structure results.

But it may well be
that we should consider also different questions.
Just as not  all group theory consist of generalization of
the Abelian case, 
so also
 there are other  natural families 
(extending the simple case), so-called
NSOP$ _2$  and NSOP$ _3 $,
on which we know basically nothing.

Probably a good test problem here is the Keisler order;
(an
e.c. $ K _1 $ is said to
be
smaller than
$ K_2 $ when for every so-called
regular  ultrafilter $ D $ in a set $ I $,
if $ M_2 ^I/D_1 $ is $ \vert I \vert^+ $-saturated
for every $ M_2 \in K_2 $ \underline{then}
this holds for $ K_1 $.)

We have a reasonable
understanding
of this order
for stable elementary classes, and we also know
that being like the theory of linear order implies maximality.
The challenge is to understand the order for simple
elementary classes and for a wider family, so-called NSOP$_3$;
we hope that this will
shed light on
those families, and lead us to a deep internal theory.

\subsection{Unstable elementary classes: dependent theories}

Simple theories include random graphs but not linear orders. On the other side
we find {\it dependent theories,}
for which the class of dense
linear orders serves as a prototype
(and dependent theories include many classes of fields, including many fields of
formal power series).

Particularly in the last decade
there has been much work
on these classes, but usually  
in more restricted contexts.

We can count the number  of so-called complete  types over $ M_1 $;
 which can be defined by:  $ a, b  \in M_2 $
realize the same type over $ M_1 $ 
where $ M_ 2 $ is an elementary
extension of $ M_1 $ if in some 
elementary extension $ N $ of $ M_2 $
 there is a automorphism  $ f $  of $ N $
over $ M_1 $ (that is , 
$ f \upharpoonright M_1 $ is the identity)
mapping $ a $ to $ b $.  
Now the class is stable when for many
cardinals $ \lambda $, if 
$ M_1 $ has $ \le \lambda $ elements, then
the number of those types is 
$ \le \lambda $, and this fails for
unstable $ T $. However,
we may %
count the above types only up to
conjugacy, that is demanding  
only that $ f $ maps $ M_1 $ onto $ M_1 $.
Now this number may be large because $ M_1 $ has too few
automorphisms, so (ignoring some points)
 we should restrict
ourselves to $ M_1 $ with enough automorphisms,
so-called saturated models.
\relax From this perspective,
for stable $ K $ , the number
is bounded (i.e., does not depend on the cardinal);
for  dependent $ K $, we get not too many;
and for  independent $ K $ we
get almost always the maximal values
$ 2^\lambda $).

A great challenge is to now understand those types, and hence dependent
 classes.

\subsection{Back to the stable setting}

There are great challenges which remain for the stable case.
\underline{The main gap conjecture
} for a family of
classes, says that for a class $ K $
(from the family),
the function $ I( \lambda , K ) $ either
is usually maximal (i.e. $ 2 ^ \lambda $)
or is  not too large , and that there is
a clear characterization.
We hope that when this is not maximal every
model can be represented by a graph
as %
 a "base" which is a tree,
with a root and the nodes are coloured by not too
many colours. 
More specifically, 
every model
can be described by 
such %
a tree of small models
put together
in a \text{``}free\text{''} %
 (hence unique) way,
the model is so called prime over this tree 
of models, but it is not claimed that the tree 
is unique.
  For general  elementary classes $ K $
we still do not know it; but if
the vocabulary is countable - we know.
Also, even for countable vocabulary, for
$ {\aleph_1} $-saturated models we do not know.
\end{document}